\documentclass[a4paper, 12pt, leqno]{amsart}
\usepackage{amsmath, mathtools}
\usepackage{amscd}
\usepackage{amssymb}
\usepackage{amsthm}
\usepackage{ascmac, color}
\usepackage{graphicx, xcolor, pdfpages}
\usepackage{mathrsfs}
\usepackage[all]{xy}
\usepackage{hyperref}

\newtheorem{lemma}{{\sc Lemma}}[section]
\newtheorem{corollary}[lemma]{{\sc Corollary}}
\newtheorem{proposition}[lemma]{{\sc Proposition}}
\newtheorem{theorem}[lemma]{{\sc Theorem}}
\theoremstyle{definition}
\newtheorem{remark}[lemma]{{\sc Remark}}

\numberwithin{equation}{section}

\def\Ga{{\mathfrak{a}}}
\def\Gb{{\mathfrak{b}}}

\def\Gg{{\mathfrak{g}}}
\def\Gh{{\mathfrak{h}}}

\def\Gl{{\mathfrak{l}}}
\def\Gm{{\mathfrak{m}}}
\def\Gn{{\mathfrak{n}}}

\def\Gp{{\mathfrak{p}}}
\def\Gr{{\mathfrak{r}}}

\def\BC{{\mathbf{C}}}

\def\BZ{{\mathbf{Z}}}

\def\CL{{\mathcal L}}

\def\CR{{\mathcal R}}

\def\CV{{\mathcal V}}

\def\SI{{\mathscr I}}

\def\ad{{\mathop{\rm ad}\nolimits}}

\def\an{{\mathop{\rm an}\nolimits}}
\def\co{\mathop{\rm co}\nolimits}

\def\deru{\partial}

\def\Ext{{\mathop{\rm Ext}\nolimits}}

\def\Hom{\mathop{\rm Hom}\nolimits}

\def\id{\mathop{\rm id}\nolimits}

\def\Ind{{\mathop{\rm Ind}\nolimits}}

\def\Image{\mathop{\rm Im}\nolimits}

\def\max{{\mathop{\rm max}\nolimits}}
\def\Mod{\mathop{\rm Mod}\nolimits}

\def\m{\Gm_{J,q}}

\def\K{{U_q(\Gl_J)}}
\def\bU{{\overline{U}}}
\def\O{{O}}
\def\bO{{\overline{O}}}
\def\F{F}
\def\al{\blacktriangleleft}
\def\ar{\blacktriangleright}

\begin{document}
\title[Differential calculus on quantized irreducible flag manifolds]
{
Differential calculus on quantized irreducible flag manifolds}
\author{Toshiyuki TANISAKI}
%\address{}
%\email{}
\begin{abstract}
We give a new description of the differential calculus on the quantized irreducible flag manifold  given by Heckenberger-Kolb.
We also define a quantum analogue of the Dolbeault complex associated to an equivariant vector bundle on the irreducible flag manifold, and show that its cohomology groups coincide with the cohomology groups given by Andersen-Polo-Wen defined using the induction functor.
\end{abstract}
\maketitle
%\tableofcontents

\section{Introduction}
\subsection{}
For a smooth complex algebraic variety $Y$ let $Y_\an$ be the associated complex manifold.
For an algebraic  vector bundle $\CL$ on $Y$ let $\CL_\an$ be the corresponding holomorphic vector bundle on $Y_\an$.
Assume that $Y$ is projective.
By Serre's GAGA \cite{Se} the cohomology group $H^r(Y;\CL)$ is naturally isomorphic to the Dolbeault cohomology group $H^{0,r}(Y_\an;\CL_\an)$.
It is an interesting theme to pursuit non-commutative analogues of this fact  (see 
\cite{BeS}).
In this paper we are concerned with the special case of this problem related to quantum groups.

Let $G$ be a connected, simply-connected simple algebraic group over the complex number field $\BC$.
Let $P$ be a parabolic subgroup of $G$, and let $L$ be a Levi subgroup of $P$.
The homogeneous space $X=G/P$  is a smooth projective variety called a (generalized) flag manifold for $G$.
For a finite-dimensional $P$-module $V$ there is associated an algebraic $G$-equivariant vector bundle $\CL_V$ on $X$, and 
we have an isomorphism 
\begin{equation}
\label{eq:isom}
H^r(X;\CL_V)\cong H^{0,r}(X_\an;\CL_{V,\an})
\end{equation}
of $G$-modules.

\subsection{}
We consider the quantum analogue of \eqref{eq:isom} using the quantum groups.
In this paper we assume that  the parameter $q$ for the quantum group is not a root of 1.
Let $O_q(G)$, $O_q(P)$ and $O_q(L)$ be the quantized coordinate algebras of $G$, $P$ and $L$ respectively.
They are Hopf algebras over $\BC$.
A quantum analogue of $H^r(X;\CL_V)$ is naturally defined using the 
derived induction functor 
\[
R^r\Ind:\{\text{$O_q(P)$-comodules}\}\to
\{\text{$O_q(G)$-comodules}\}
\]
(see \cite{APW}).
In fact $R^r\Ind(V_q)$ for an $O_q(P)$-comodule $V_q$ plays the role of $H^r(X;\CL_V)$ as above.
However, it is not so clear what is the counterpart of $H^{0,r}(X_\an;\CL_{V,\an})$ in the quantum setting.
Take a maximal compact subgroup $G_u$ of $G_\an$ such that $L_u=G_u\cap P_\an$ is a maximal compact subgroup of  $L_\an$.
As a real manifold we have $X_\an=G_u/L_u$.
Hence the subalgebra $O_q(G/L)$ of $O_q(G)$ consisting of $O_q(L)$-covariant vectors is regarded as a quantum analogue of the $G_u$-finite part of the algebra of 
$\BC$-valued $C^\infty$-functions on $X_{\an}$.

Assume that the unipotent radical of $P$ is commutative 
(in this case $G/P$ is called an irreducible flag manifold in some references).
In \cite{HK1}, \cite{HK2}, \cite{HK3} Heckenberger-Kolb  introduced a quantum analogue of (the $G_u$-finite part of) the Dolbeault complex
\begin{equation}
\label{eq:Dol1}
\Omega^{0,\bullet}_{X_\an}=
[0\to \Omega^{0,0}_{X_\an}\xrightarrow{\overline{\deru}}\Omega_{X_\an}^{0,1}\xrightarrow{\overline{\deru}}\cdots].
\end{equation}
We can regard $O_q(G/L)$ equipped with this analogue of \eqref{eq:Dol1} as a quantum analogue of the complex manifold $X_\an$.
Moreover, Carotenutto-D\'{i}az Garc\'{i}a-\'{O} Buachalla \cite{CDO} defined, for a one-dimensional $O_q(L)$-comodule $V_q$, a quantum analogue of (the $G_u$-finite part of) the Dolbeault complex
\begin{equation}
\label{eq:Dol2}
\Omega^{0,\bullet}_{X_\an} \otimes_{O_{X_\an}} \CL_{V,\an}=
[0\to 
\Omega^{0,0}_{X_\an}\otimes_{O_{X_\an}}  \CL_{V,\an}
\xrightarrow{\nabla}
\Omega_{X_\an}^{0,1}\otimes_{O_{X_\an}}  \CL_{V,\an}
\xrightarrow{\nabla}\cdots]
\end{equation}
associated to a line bundle $\CL_{V,\an}$, and proved that the Borel-Weil type theorem holds for the quantum analogue of $H^{0,0}(X_\an;\CL_{V,\an})$ (in the case $V_q$ is one-dimensional).

\subsection{}
In the present paper still assuming that the unipotent radical of $P$ is commutative, we 
give a new description of the quantum analogue of \eqref{eq:Dol1}.
We also show that a quantum analogue of \eqref{eq:Dol2} is defined 
for any $O_q(P)$-comodule not necessarily one-dimensional.
Finally, we show that the quantum analogue of \eqref{eq:isom} holds for any $O_q(P)$-comodule $V_q$.
In particular,
the Borel-Weil-Bott type theorem holds for the quantum analogue of $H^{0,r}(X_\an;\CL_{V,\an})$ for any irreducible $O_q(P)$-comodule $V_q$.

\subsection{}
In this paper the notation $\otimes$ without a subscript stands for $\otimes_\BC$.

Let $H$ be a Hopf algebra over $\BC$.
We denote the comultiplication, the counit and the antipode of $H$ by
\[
\Delta:H\to H\otimes H, 
\qquad
\varepsilon:H\to\BC,
\qquad
S:H\to H
\]
respectively.
The category of left (resp.\ right) $H$-comodules
is denoted by ${}^H\!\Mod$ (resp.\ $\Mod^H$).
For $h\in H$ we set $h^+=h-\varepsilon(h)1$.
For a subalgebra $K$ of $H$ we set
\[
K^+=\{k\in K\mid \varepsilon(k)=0\}.
\]
We use Sweedler's notation 
\[
\Delta^{(r)}(h)=\sum h_{(0)}\otimes\cdots\otimes h_{(r)}
\]
for the iterated comultiplication $\Delta^{(r)}:H\to H^{\otimes r+1}$.
Similarly, for $N\in{}^H\!\Mod$ 
(resp.\ $M\in\Mod^H$) we write
\[
N\to H^{\otimes r}\otimes N\quad(n\mapsto\sum n_{(-r)}\otimes\dots\otimes n_{(-1)}\otimes n_{(0)})
\]
\[
\left(\text{resp.\ }
M\to M\otimes H^{\otimes r}\quad(m\mapsto\sum m_{(0)}\otimes m_{(1)}\otimes\cdots\otimes m_{(r)})
\right)
\]
for the iterated coactions.
For $M\in\Mod^H$ we set
\[
M^{\co(H)}=
\{m\in M\mid
\sum m_{0}\otimes m_{(1)}=m\otimes 1\}.
\]
For $M\in\Mod^H$ and $N\in{}^H\!\Mod$  we denote by $M\Box_HN$ the cotensor product of $M$ and $N$ defined by 
\begin{align*}
&M\Box_HN
\\
=&\{\sum m_r\otimes n_r\in M\otimes N\mid
\sum m_{r(0)}\otimes m_{r(1)}\otimes n_r
=\sum m_r\otimes n_{r(-1)}\otimes n_{r(0)}
\}.
\end{align*}
We have
\[
M\Box_HN=(M\otimes N)^{\co(H)}
\]
with respect to the right $H$-comodule structure of $M\otimes N$ given by
\[
M\otimes N\to (M\otimes N)\otimes H
\qquad(m\otimes n\mapsto (m_{(0)}\otimes n_{(0)})
\otimes m_{(1)}(Sn_{(-1)})).
\]

\section{Quantum groups}
\subsection{}
Let $G$ be a connected, simply-connected simple algebraic group over $\BC$, and let $H$ be a maximal torus of $G$.
We denote by $\Gg$ and $\Gh$ the Lie algebras of $G$ and $H$ respectively.
Let $\Delta\subset\Gh^*$ be the set of roots for $(\Gg,\Gh)$.
We fix a set of simple roots $\{\alpha_i\mid i\in I\}$, and denote by $\Delta^+$ the corresponding set of positive roots.
We denote by $W\subset GL(\Gh^*)$ the Weyl group.
It is generated by simple reflections $s_i$ ($i\in I$).
We denote by $w_0$ the longest element of $W$.
Let 
\begin{equation}
\label{eq:bil}
(\;,\;):\Gh^*\times\Gh^*\to\BC
\end{equation}
be the $W$-invariant symmetric bilinear form such that $(\alpha,\alpha)=2$ for short roots.
We set $a_{ij}=2(\alpha_i,\alpha_j)/(\alpha_i,\alpha_i)$ for $i, j\in I$.
%For $\alpha\in\Delta$ the corresponding coroot $\alpha^\vee$ is defined by $\alpha^\vee=2\alpha/(\alpha,\alpha)^{-1}$.
Let $Q$ 
%and $Q^\vee$ 
be the root lattice 
%and the coroot lattice respectively, 
and let $\Lambda\subset\Gh^*$ be the set of integral weights.
We have $Q\subset \Lambda$.
Moreover, 
by our choice of \eqref{eq:bil}
we have 
%$Q\subset Q^\vee$ and 
$(\Lambda, Q)\subset\BZ$.
We set
\[
Q^+=\sum_{i\in I}\BZ_{\geqq0}\alpha_i.
%\qquad
%\Lambda^+=
%\{\lambda\in\Lambda\mid
%(\lambda,\alpha^\vee)\geqq0\;\text{ for } \alpha\in\Delta^+\}.
\]
For $\alpha\in\Delta$ we denote by $\Gg_\alpha$ the root subspace.
We define sublagebras $\Gn$, $\Gb$ of $\Gg$ by
\[
\Gn=\bigoplus_{\alpha\in\Delta^+}\Gg_\alpha,
\qquad
\Gb=\Gh\oplus\Gn.
\]
%We define sublagebras $\Gn$, $\Gn^-$ $\Gb$, $\Gb^-$ of $\Gg$ by
%\[
%\Gn=\bigoplus_{\alpha\in\Delta^+}\Gg_\alpha,
%\qquad
%\Gn^-=\bigoplus_{\alpha\in\Delta^+}\Gg_{-\alpha},
%\qquad
%\Gb=\Gh\oplus\Gn,
%\qquad
%\Gb^-=\Gh\oplus\Gn^-.
%\]

Let  $J$ be a subset of $I$. 
We set
\[
\Delta_J=\Delta\cap\left(
\sum_{j\in J}\BZ\alpha_j\right),
\qquad
\Delta_J^+=\Delta^+\cap\Delta_J,
\qquad
\Xi_J=\Delta^+\setminus\Delta_J^+.
\]
We denote the longest element of $W_J=\langle s_j\mid j\in J\rangle$ by $w_J$.
We define subalgebras $\Gl_J$, $\Gb_J$, $\Gm_J$, $\Gp_J$ of $\Gg$ by
\[
\Gl_J=\Gh\oplus\left(\bigoplus_{\alpha\in\Delta_J}\Gg_\alpha\right),
\quad
\Gb_J=\Gh\oplus\left(\bigoplus_{\alpha\in\Delta^+_J}\Gg_\alpha\right),
\quad
\Gm_J=\bigoplus_{\alpha\in\Xi_J}\Gg_\alpha,
\quad
\Gp_J=\Gl_J\oplus\Gm_J.
\]
We denote by $L_J$, $P_J$ the closed connected subgroups of $G$ with Lie algebras $\Gl_J$, $\Gp_J$ respectively.
Let $\sum_{i\in I}m_i\alpha_i$ be the highest root.
The following fact is well-known.
%%%%%%%%%%%%%%%%%%
\begin{lemma}
\label{lem:irred}
The following conditions for a subset $J$ of $I$ are equivalent to each other.
\begin{itemize}
\item[(a)] 
$J\ne I$ and $[\Gm_J,\Gm_J]=\{0\}$.
\item[(b)] 
$\Gm_J$ is an irreducible $\Gl_J$-module with respect to the adjoint action.
\item[(c)] 
$J=I\setminus\{i_0\}$ with $m_{i_0}=1$.

\end{itemize}
\end{lemma}
%%%%%%%%%%%%%%%%%%

%%%%%%%%%%%%%%%%%%%
\subsection{}
We fix $q\in\BC^\times$.
In this paper we assume that $q$ is not a root of $1$.
We denote by $U_q(\Gg)$ the quantized enveloping algebra of $\Gg$.
It is a $\BC$-algebra generated by the elements 
\[
e_i, \; f_i\quad(i\in I),\qquad k_\gamma\quad(\gamma\in Q)
\]
satisfying the relations:
\begin{align*}
&k_0=1,\qquad
k_\gamma k_{\gamma'}=k_{\gamma+\gamma'}
\qquad&(\gamma, \gamma'\in Q),
\\
&k_\gamma e_i=q^{(\gamma,\alpha_i)}e_i k_\gamma,
\quad
k_\gamma f_i=q^{-(\gamma,\alpha_i)}f_i k_\gamma,
\qquad
&(\gamma\in Q,\; i\in I),
\\
&e_if_j-f_je_i=\delta_{ij}\frac{k_i-k_i^{-1}}{q_i-q_i^{-1}}
\qquad&(i, j\in I),
\\
&
\sum_{r=0}^{1-{a_{ij}}}
(-1)^{r}e_i^{(1-a_{ij}-r)}e_je_i^{(r)}=
\sum_{r=0}^{1-{a_{ij}}}
(-1)^{r}f_i^{(1-a_{ij}-r)}f_jf_i^{(r)}=0
\qquad&(i, j\in I,\; i\ne j),
\end{align*}
where
$k_i=k_{\alpha_i}$, $q_i=q^{(\alpha_i,\alpha_i)/2}$ and
\[
e_i^{(r)}=
\left(\prod_{s=1}^r\frac{q_i^s-q_i^{-s}}{q_i-q_i^{-1}}\right)^{-1}
e_i^r,
\qquad
f_i^{(r)}=
\left(\prod_{s=1}^r\frac{q_i^s-q_i^{-s}}{q_i-q_i^{-1}}\right)^{-1}
f_i^r
\]
for $i\in I$ and $r\geqq0$.
Its Hopf algebra structure is given by 
\[
\Delta(e_i)=e_i\otimes k_i+1\otimes e_i,
\quad
\Delta(f_i)=f_i\otimes 1+k_i^{-1}\otimes f_i,
\quad
\Delta(k_\gamma)=k_\gamma\otimes k_\gamma,
\]
\[
\varepsilon(e_i)=\varepsilon(f_i)=0,\quad\varepsilon(k_\gamma)=1,
\]
\[
S(e_i)=-e_ik_i^{-1},\quad
S(f_i)=-k_if_i,\quad
S(k_\gamma)=k_{-\gamma}
\]
for $i\in I$, $\gamma\in Q$.
We define subalgebras $U_q(\Gh)$, $U_q(\Gn)$, 
$U_q(\Gb)$  of $U_q(\Gg)$ by
\[
U_q(\Gh)=\langle k_\gamma\mid \gamma\in Q\rangle,
\qquad
U_q(\Gn)=\langle e_i\mid i\in I\rangle,
\qquad
U_q(\Gb)=\langle U_q(\Gh), U_q(\Gn)\rangle.
\]
Then $U_q(\Gh)$, $U_q(\Gb)$ are  Hopf subalgebras of $U_q(\Gg)$.
For $\lambda\in\Lambda$ we define a character $\chi_\lambda:U_q(\Gh)\to \BC$ by 
\[
\chi_\lambda(k_\gamma)=q^{(\lambda,\gamma)}
\qquad(\gamma\in Q).
\]
For a left (resp.\ right) $U_q(\Gh)$-module $M$ and $\lambda\in\Lambda$ we set 
\[
M_\lambda=\{m\in M\mid hm=\chi_\lambda(h)m\;
(\text{resp.} \;
mh=\chi_\lambda(h)m)\;\;
\text{for } h\in U_q(\Gh)\}.
\]
We say that a left or right $U_q(\Gh)$-module is a weight module if we have the weight space decomposition
\[
M=\bigoplus_{\lambda\in\Lambda}M_\lambda.
\]

We define the adjoint action of $U_q(\Gg)$ on $U_q(\Gg)$ by
\[
\ad(u)(x)=\sum u_{(1)}x(S^{-1}u_{(0)})
\qquad (u, x\in U_q(\Gg)).
\]
The subalgebra $U_q(\Gn)$ is a left $U_q(\Gh)$-module with respect to the adjoint action with the weight space decomposition 
\[
U_q(\Gn)=\bigoplus_{\gamma\in Q^+}U_q(\Gn)_\gamma.
\]

Following Lusztig \cite{Lb} we define an algebra automorphism $T_i$ ($i\in I$) of $U_q(\Gg)$ by
\begin{align*}
T_i(e_j)=&
\begin{cases}
-f_ik_i\quad &(i=j)
\\
\sum_{r=0}^{-a_{ij}}(-q_i)^{-r}e_i^{(-a_{ij}-r)}e_je_i^{(r)}
\quad &(i\ne j),
\end{cases}
\\
T_i(f_j)=&
\begin{cases}
-k_i^{-1}e_i\quad &(i=j)
\\
\sum_{r=0}^{-a_{ij}}(-q_i)^{r}f_i^{(r)}f_jf_i^{(-a_{ij}-r)}
\quad &(i\ne j),
\end{cases}
\\
T_i(k_\gamma)=&k_{s_i\gamma}
\end{align*}
for $j\in I$, $\gamma\in Q$.
Its inverse is given by
\begin{align*}
T_i^{-1}(e_j)=&
\begin{cases}
-k_i^{-1}f_i\quad &(i=j)
\\
\sum_{r=0}^{-a_{ij}}(-q_i)^{-r}e_i^{(r)}e_je_i^{(-a_{ij}-r)}
\quad &(i\ne j),
\end{cases}
\\
T_i^{-1}(f_j)=&
\begin{cases}
-k_i^{-1}e_i\quad &(i=j)
\\
\sum_{r=0}^{-a_{ij}}(-q_i)^{r}f_i^{(-a_{ij}-r)}f_jf_i^{(r)}
\quad &(i\ne j),
\end{cases}
\\
T_i^{-1}(k_\gamma)=&k_{s_i\gamma}.
\end{align*}
For $w\in W$ choose a reduced expression 
$w=s_{i_1}\dots s_{i_r}$
of $w$ and define an algebra automorphism $T_w$ of $U_q(\Gg)$ by 
$
T_w=T_{i_1}\dots T_{i_r}.
$
It does not depend on the choice of a reduced expression.

\subsection{}
Let $J\subset I$.
We define subalgebras $U_q(\Gl_J)$, $U_q(\Gb_J)$, $U_q(\Gp_J)$, $U_q(\Gm_J)$ of $U_q(\Gg)$ by
\[
U_q(\Gl_J)=
\langle
U_q(\Gh), e_j, f_j\mid j\in J\rangle,
\qquad
U_q(\Gb_J)=
\langle
U_q(\Gh), e_j\mid j\in J\rangle,
\]
\[
U_q(\Gp_J)=
\langle
U_q(\Gh), e_i, f_j\mid i\in I, j\in J\rangle,
\qquad
U_q(\Gm_J)=U_q(\Gn)\cap T_{w_J}(U_q(\Gn)).
\]
Then $U_q(\Gl_J)$, $U_q(\Gb_J)$ and $U_q(\Gp_J)$ are Hopf subalgebras of $U_q(\Gg)$.

We have the following (see for example \cite[Proposition 2.10]{TM}).
%%%%%%%%%%%%%%%%%%
\begin{lemma}
\label{lem:dec}
The multiplication of $U_q(\Gb)$ induces isomorphisms
\[
U_q(\Gm_J)\otimes U_q(\Gb_J)\cong U_q(\Gb)
\cong U_q(\Gb_J)\otimes U_q(\Gm_J)
\]
of vector spaces.
\end{lemma}
%%%%%%%%%%%%%%%%%%

We have also the following (see for example \cite[Lemma 2.8]{TM}).
%%%%%%%%%%%%%%%%%%
\begin{lemma}
\label{lem:DUm}
$\Delta(U_q(\Gm_J))\subset U_q(\Gm_J)\otimes U_q(\Gb)$.
\end{lemma}
%%%%%%%%%%%%%%%%%%

Similarly to \cite{KMT} 
we have the following.
%%%%%%%%%%%%%%%%%%
\begin{lemma}
\label{lem:adUm}
$\ad(U_q(\Gl_J))(U_q(\Gm_J))\subset U_q(\Gm_J)$.
\end{lemma}

We choose a reduced expression 
$w_Jw_0=s_{i_1}\dots s_{i_d}$ of $w_Jw_0$ and set
\begin{equation}
\label{eq:beta}
\beta_r=s_{i_1}\dots s_{i_{r-1}}(\alpha_{i_r})
\qquad(r=1,\dots d).
\end{equation}
Then we have $\Xi_J=\{\beta_1,\dots, \beta_d\}$.
Set 
\begin{equation}
\label{eq:b}
e_{\beta_r}=T_{i_1}^{-1}\dots T_{i_{r-1}}^{-1}(e_{i_r})
\qquad(r=1,\dots d).
\end{equation}
Then we have the following (see \cite[Proposition 40.2.1]{Lb}).
\begin{lemma}
\label{lem:basis}
The set
\[
\{
e_{\beta_1}^{n_1}\dots e_{\beta_d}^{n_d}\mid
(n_1,\dots,n_d)\in(\BZ_{\geqq0})^d\}
\]
forms a basis of $U_q(\Gm_J)$.
\end{lemma}
%%%%%%%%%%%%%%%%%%%%
\subsection{}
For $\Ga=\Gg$, $\Gh$, $\Gp_J$, $\Gl_J$ we say that a left (resp.\ right) $U_q(\Ga)$-module $M$ is integrable if it is a weight module as a $U_q(\Gh)$-module, 
and if we have $\dim U_q(\Ga)m<\infty$
(resp.\ $\dim mU_q(\Ga)<\infty$)
for any $m\in M$ .
Let $A$ be the closed connected subgroup of $G$ with Lie algebra $\Ga$.
We denote by $O_q(A)$ the subspace of $U_q(\Ga)^*$ spanned by the matrix coefficients of  integrable left $U_q(\Ga)$-modules 
(it coincides with the subspace of $U_q(\Ga)^*$  spanned by the matrix coefficients of  integrable right $U_q(\Ga)$-modules).
Then $O_q(A)$ turns out to be a Hopf algebra whose multiplication, comultiplication, unit, counit, antipode are given by the transpose of the 
comultiplication, multiplication, counit, unit, antipode of $U_q(\Ga)$ respectively.
The canonical pairing
\begin{equation}
\label{eq:dual}
\langle\;,\;\rangle:O_q(A)\times U_q(\Ga)\to\BC
\end{equation}
is non-degenerate in the sense that 
$u\in U_q(\Ga)$ satisfying $\langle O_q(A),u\rangle=\{0\}$ is zero, and
$\varphi\in O_q(A)$ satisfying $\langle \varphi, U_q(\Ga)\rangle=\{0\}$ is zero.
We have a $U_q(\Ga)$-bimodule structure of $O_q(A)$ given by
\[
\langle u_1\varphi u_2,u\rangle=
\langle\varphi,u_2uu_1\rangle
\qquad(\varphi\in O_q(A), u_1, u_2, u\in U_q(\Ga)).
\]
For $u\in U_q(\Ga)$, $\varphi\in O_q(A)$ we have
\begin{equation}
\label{eq:Delta1}
\Delta(u\varphi)=\sum\varphi_{(0)}\otimes u\varphi_{(1)},
\qquad
\Delta(\varphi u)=\sum\varphi_{(0)}u\otimes \varphi_{(1)},
\end{equation}
\begin{equation}
\label{eq:Delta2}
\sum u\varphi_{(0)}\otimes \varphi_{(1)}
=\sum\varphi_{(0)}\otimes \varphi_{(1)}u.
\end{equation}

The category of integrable left (resp.\ right) $U_q(\Ga)$-modules is naturally equivalent to the category of 
right (resp.\ left) $O_q(A)$-comodules $\Mod^{O_q(A)}$ (resp.\ ${}^{O_q(A)}\!\Mod$). 

For $J\subset I$ we have a sequence 
\[
O_q(G)\twoheadrightarrow
O_q(P_J)\twoheadrightarrow
O_q(L_J)
\]
of surjective Hopf algebra homomorphisms.

%%%%%%%%%%%%%%%%%%%%%%%
\section{Differential calculus over $O_q(G)$}
\label{sec:O}
\subsection{}
In the rest of this paper we fix a subset $J$ of $I$ satisfying the equivalent conditions in Lemma \ref{lem:irred}.
In particular, we have $J=I\setminus\{i_0\}$ for some $i_0\in I$.
Define a numbering $\Xi_J=\{\beta_1,\dots, \beta_d\}$ of elements of $\Xi_J$ as in \eqref{eq:beta}.
For each $\beta\in\Xi_J$ we have $e_\beta\in U_q(\Gm_J)$ given  by \eqref{eq:b}.
We set 
\begin{equation}
\m=\sum_{\beta\in\Xi_J}\BC e_\beta
\subset U_q(\Gm_J).
\end{equation}

Consider the grading 
\begin{equation}
\label{eq:Umgr}
U_q(\Gm_J)=\bigoplus_{n=0}^\infty U_q(\Gm_J)(n)
\end{equation}
of the algebra $U_q(\Gm_J)$, where $U_q(\Gm_J)(n)$
is the sum of the weight spaces whose weight belongs  to 
$n\alpha_{i_0}+\sum_{j\in J}\BZ_{\geqq0}\alpha_j$.
By the condition (c) in Lemma \ref{lem:irred} and Lemma \ref{lem:basis}
we have 
$U_q(\Gm_J)(0)=\BC$,
$U_q(\Gm_J)(1)=\m$, and $U_q(\Gm_J)(n)$ coincides with the image of 
$\m^{\otimes n}\to U_q(\Gm_J)$ defined by the multiplication of $U_q(\Gm_J)$.
Considering the weights we obtain from Lemma \ref{lem:DUm} the following.
%%%%%%%%%%%%%%%%%%

\begin{lemma}
\label{lem:Dm}
We have
\[
\Delta(x)\in 1\otimes x+\m\otimes U_q(\Gb_J)
\qquad
(x\in\m).
\]
In particular, we have
\[
\sum x_{(0)}^+\otimes x_{(1)}\in \m\otimes \K\qquad
(x\in \m).
\]
\end{lemma}

%%%%%%%%%%%%%%%%%%
\begin{corollary}
\label{cor:Dm}
For $x\in\m$ we have
\[
\Delta(x)-(x\otimes1+1\otimes x)\in\m\otimes U_q(\Gb_J)^+,
\qquad
Sx+x\in U_q(\Gb)U_q(\Gb_J)^+.
\]
\end{corollary}
%%%%%%%%%%%%%%%%%%
\begin{proof}
By Lemma \ref{lem:Dm} and the definition of the comultiplication of $U_q(\Gg)$ we have
\begin{align*}
\Delta(e_\beta)\in &
1\otimes e_\beta+e_\beta\otimes k_\beta
+\m\otimes U_q(\Gb_J)^+
\\
=&1\otimes e_\beta+e_\beta\otimes 1
+e_\beta\otimes(k_\beta-1)
+\m\otimes U_q(\Gb_J)^+
\\
=&1\otimes e_\beta+e_\beta\otimes 1
+\m\otimes U_q(\Gb_J)^+,
\end{align*}
and hence the first formula holds.
The second one follows by applying $S\otimes\id$
to the first.
\end{proof}
We also obtain from Lemma \ref{lem:adUm} the following.
\begin{lemma}
\label{lem:adm}
$\ad(\K)(\m)\subset\m$.
\end{lemma}

By Lemma \ref{lem:basis} $U_q(\Gm_J)(2)$ has a basis,
\[
\{e_{\beta_m}e_{\beta_n}\mid 1\leqq m\leqq n\leqq d\}.
\]
Hence for $d\geqq r>s\geqq 1$ we can write 
\begin{equation}
\label{eq:rel}
e_{\beta_r}e_{\beta_s}=
\sum_{1\leqq m\leqq n\leqq d}C^{r,s}_{m,n}e_{\beta_m}e_{\beta_n}
\end{equation}
for some $C^{r,s}_{m,n}\in\BC$.
By Lemma \ref{lem:basis} we see that the algebra $U_q(\Gm_J)$ is generated by the elements $e_\beta$ ($\beta\in\Xi_J$) satisfying the defining equations \eqref{eq:rel} for $d\geqq r>s\geqq 1$.
In particular, $U_q(\Gm_J)$ is a quadratic algebra generated by $\m$.
We set
\begin{align}
\Gr=&
\{\sum_c x_c\otimes y_c\in\m\otimes\m\mid
\sum x_cy_c=0\}
\\
\nonumber
=&\sum_{r>s}\BC(e_{\beta_r}\otimes e_{\beta_s}-
\sum_{m\leqq n}C^{r,s}_{m,n}e_{\beta_m}\otimes e_{\beta_n})
\subset \m\otimes\m.
\end{align}

We denote by $\wedge\m^*$ the quadratic algebra dual to $U_q(\Gm_J)$.
Namely, we set
\begin{equation}
\label{eq:Wm}
\wedge\m^*=T\m^*/(T\m^*)\Gr^\perp (T\m^*),
\end{equation}
where $T\m^*$ is the tensor algebra of the $\BC$-module $\m^*$ and 
\[
\Gr^\perp=\{
z\in\m^*\otimes\m^*\mid \langle z,\Gr\rangle=\{0\}\}.
\]
The multiplication of $\wedge\m^*$ is written as 
\[
\wedge\m^*\times \wedge\m^*\to \wedge\m^*
\qquad((a,b)\mapsto a\wedge b).
\]
For $n\geqq0$ we denote by $\wedge^n\m^*$ the image of $(\m^*)^{\otimes n}\to \wedge\m^*$ defined by the multiplication of $\wedge\m^*$.
Then we have 
\[
\wedge\m^*=\bigoplus_n\wedge^n\m^*.
\]
Denote by $\{e^*_{\beta_1},\dots, e^*_{\beta_d}\}$ the basis of $\m^*$ dual to $\{e_{\beta_1},\dots, e_{\beta_d}\}$.
By a general fact on quadratic algebras (see \cite[Theorem 4.1]{PP}) we have 
$\wedge^n\m^*=\{0\}$ for $n>d(d-1)/2$, and 
$\wedge^n\m^*$ for $0\leqq n\leqq d(d-1)/2$ has a basis
\[
\{e^*_{\beta_{r_1}}\wedge\dots\wedge e^*_{\beta_{r_n}}\mid
d\geqq r_1>\dots>r_n\geqq 1\}.
\]
%%%%%%%%%%%%%%%%%%
\begin{lemma}
\label{lem:eei}
We have 
\[
\sum_{\beta,\gamma\in\Xi_J}e_\beta e_\gamma\otimes (e^*_\beta\wedge e^*_\gamma)=0
\]
in $U_q(\Gm)\otimes \wedge\m$.
\end{lemma}
%%%%%%%%%%%%%%%%%%
\begin{proof}
For $\sum_c x_c\otimes y_c\in\Gr$
we have
\[
\sum_c\sum_{\beta,\gamma\in\Xi_J} \langle e^*_\beta\otimes e^*_\gamma, x_c\otimes y_c\rangle e_\beta e_\gamma
=\sum_c(\sum_\beta\langle e^*_\beta,x_c\rangle e_\beta)
(\sum_\gamma\langle e^*_\gamma,y_c\rangle e_\gamma)
=\sum_cx_cy_c=0,
\]
and hence 
$\sum_{\beta,\gamma\in\Xi_J}e_\beta e_\gamma\otimes (e^*_\beta\otimes e^*_\gamma)\in U_q(\Gm)\otimes\Gr^\perp$.
Therefore, $\sum_{\beta,\gamma\in\Xi_J}e_\beta e_\gamma\otimes (e^*_\beta\wedge e^*_\gamma)=0$
in $U_q(\Gm)\otimes \wedge\m$.
\end{proof}
\begin{proposition}
\label{prop:identity}
For $\sum_cx_c\otimes y_c\in\Gr$ we have
\begin{align}
\label{eq:Id1}
&\sum x^+_{c(0)}y^+_{c(0)}\otimes x_{c(1)}y_{c(1)}=0,
\\
\label{eq:Id2}
&\sum x^+_{c(0)}\otimes x_{c(1)}y_c
+\sum y^+_{c(0)}\otimes x_cy_{c(1)}=0,
\\
\label{eq:Id3}
&\sum x_c\otimes y_c
+\sum \ad(x_{c(1)})(y_c)\otimes x_{c(0)}^+=0.
\end{align}
\end{proposition}
\begin{proof}
By $\sum x_cy_c=0$ we have
\begin{align*}
0=&\sum x_{c(0)}y_{c(0)}\otimes x_{c(1)}y_{c(1)}
\\
=&\sum(x^+_{c(0)}+\varepsilon(x_{c(0)}))
(y^+_{c(0)}+\varepsilon(y_{c(0)}))
\otimes x_{c(1)}y_{c(1)}
\\
=&
\sum x^+_{c(0)}y^+_{c(0)}\otimes x_{c(1)}y_{c(1)}
+\sum x^+_{c(0)}\otimes x_{c(1)}y_c
+\sum y^+_{c(0)}\otimes x_cy_{c(1)}
+1\otimes \sum x_cy_c
\\
=&
\sum x^+_{c(0)}y^+_{c(0)}\otimes x_{c(1)}y_{c(1)}
+\sum x^+_{c(0)}\otimes x_{c(1)}y_c
+\sum y^+_{c(0)}\otimes x_cy_{c(1)}.
\end{align*}
Hence by
\begin{align*}
\sum x^+_{c(0)}\otimes x_{c(1)}y_c
+\sum y^+_{c(0)}\otimes x_cy_{c(1)}\in&U_q(\Gm_J)(1)\otimes U_q(\Gg),
\\
\sum x^+_{c(0)}y^+_{c(0)}\otimes x_{c(1)}y_{c(1)}
\in&U_q(\Gm_J)(2)\otimes U_q(\Gg)
\end{align*}
and \eqref{eq:Umgr} we obtain \eqref{eq:Id1}, \eqref{eq:Id2}.

We have
\begin{align*}
&\sum(\ad(x_{c(2)})(y_c))x_{c(1)}\otimes x_{c(0)}^+=
\sum x_{c(3)}y_c(S^{-1}x_{c(2)})x_{c(1)}\otimes x_{c(0)}^+
\\
=&
\sum x_{c(1)}y_c\otimes x_{c(0)}^+,
\end{align*}
and hence
\[
\sum(\ad(x_{c(2)})(y_c))x_{c(1)}\otimes x_{c(0)}^+
+\sum x_cy_{c(1)}\otimes y_{c(0)}^+=0
\]
by \eqref{eq:Id2}.
Note that 
\[
x_c, \ad(x_{c(2)})(y_c)\in \m,\quad
y_{c(1)}, x_{c(1)}\in U_q(\Gb_J).
\]
Hence by Lemma \ref{lem:dec} 
we obtain
\[
\sum\ad(x_{c(2)})(y_c)\otimes x_{c(1)}\otimes x_{c(0)}^+
+\sum x_c\otimes y_{c(1)}\otimes y_{c(0)}^+=0
\]
in $\m\otimes U_q(\Gb_J)\otimes \m$.
We easily obtain \eqref{eq:Id3} by applying $\id\otimes\varepsilon\otimes\id$ to it.
\end{proof}
%%%%%%%%%%%%%%%%%%%%%%%%
\subsection{}
Since the pairing $O_q(G)\times U_q(\Gg)\to\BC$ is non-degenerate and since $\m$ is a finite-dimensional subspace of $U_q(\Gg)$, the natural linear map $O_q(G)\to\m^*$ given by the restriction is surjective.
Hence $O_q(G)^+\to\m^*$ is also surjective.
We see by Lemma \ref{lem:Dm} that the kernel $\SI_0$ of $O_q(G)^+\to\m^*$ is a right ideal of $O_q(G)$.
We regard $\m^*$ as a right $O_q(G)$-module by identifying it with $O_q(G)^+/\SI_0$.
%%%%%%%%%%%%%%%%%%

\begin{lemma}
\label{lem:mO}
The right $O_q(G)$-module structure of $\m^*$ is naturally extended to that of $\wedge\m^*$ by
\begin{itemize}
\item[(a)]
$1\varphi=\varepsilon(\varphi)1$\quad$(\varphi\in O_q(G))$,
\item[(b)]
$(a\wedge b)\varphi=\sum a\varphi_{(0)}\wedge b\varphi_{(1)}$
\quad
$(a, b\in\wedge\m^*, \varphi\in O_q(G))$.

\end{itemize}
\end{lemma}
\begin{proof}
The tensor algebra $T\m^*$ is a right $O_q(G)$-module by
\[
(a_0\cdots a_n)\varphi=\sum(a_0\varphi_{(0)})\dots(a_n\varphi_{(n)})
\qquad(a_0,\cdots, a_n\in\m^*,\; \varphi\in O_q(G)).
\]
Hence by \eqref{eq:Wm} it is sufficient to show that 
$\Gr^\perp$ is a right $O_q(G)$-submodule of 
$\m^*\otimes \m^*\subset T\m^*$.
Let 
$\sum b_r\otimes b'_r\in\Gr^\perp$ and $\varphi\in O_q(G)$.
For $\sum x_c\otimes y_c\in\Gr$ we have
\begin{align*}
&\sum \langle b_r\varphi_{(0)}\otimes b'_r\varphi_{(1)},x_c\otimes y_c\rangle
=
\sum \langle b_r\varphi_{(0)},x_c\rangle
\langle b'_r\varphi_{(1)},y_c\rangle
\\
=&
\sum \langle b_r, x_{c(0)}^+\rangle
\langle\varphi_{(0)},x_{c(1)}\rangle
\langle b'_r, y_{c(0)}^+\rangle
\langle\varphi_{(1)},y_{c(1)}\rangle
\\
=&
\sum \langle b_r\otimes b'_r\otimes\varphi, 
x_{c(0)}^+\otimes y_{c(0)}^+\otimes x_{c(1)}y_{c(1)}
\rangle=0
\end{align*}
by \eqref{eq:Id1}.
\end{proof}
%%%%%%%%%%%%%%%%%%
Hence we obtain a graded algebra structure of $O_q(G)\otimes \wedge\m^*=\bigoplus_{n=0}^{d(d-1)/2} O_q(G)\otimes \wedge^n\m^*$ given by
\begin{equation}
(\varphi\otimes a)\wedge(\psi\otimes b)
=
\sum \varphi\psi_{(0)}\otimes (a\psi_{(1)})\wedge b
\qquad(\varphi, \psi\in O_q(G), a, b\in\wedge\m^*).
\end{equation}
Note that $O_q(G)$ and $\wedge\m^*$ are naturally identified with subalgebras of $O_q(G)\otimes \wedge\m^*$ by the embeddings
$O_q(G)\cong O_q(G)\otimes 1\subset O_q(G)\otimes \wedge\m^*$ and 
$\wedge\m^*\cong1\otimes \wedge\m^*\subset O_q(G)\otimes \wedge\m^*$.
In particular, $O_q(G)\otimes \wedge\m^*$ is naturally an $O_q(G)$-bimodule by 
\begin{equation}
\psi(\varphi\otimes a)\psi'
=\sum\psi\varphi\psi'_{(0)}\otimes a\psi'_{(1)}
\qquad(\varphi, \psi, \psi'\in O_q(G), a\in\wedge\m^*).
\end{equation}

By Lemma \ref{lem:adm} $\m$ is a left $\K$-module via the adjoint action.
Hence we have a right $\K$-module structure 
\begin{equation}
\m^*\otimes\K\to\m^*
\qquad(a\otimes k\mapsto a\al k)
\end{equation}
of $\m^*$ given by
\[
\langle a\al k,x\rangle=\langle a,\ad(k)(x)\rangle
\qquad(a,\in\m^*, k\in\K, x\in\m).
\]
%%%%%%%%%%%%%%%%%%

\begin{lemma}
\label{lem:mK}
The right $\K$-module structure of $\m^*$ is naturally extended to that of $\wedge\m^*$ by
\begin{itemize}
\item[(a)]
$1\al k=\varepsilon(k)1$\quad$(k\in\K)$,
\item[(b)]
$(a\wedge b)\al k=\sum(a\al k_{(1)})\wedge(b\al k_{(0)})$
\quad
$(a, b\in\wedge\m^*, k\in\K)$.
\end{itemize}
\end{lemma}
\begin{proof}
The tensor algebra $T\m^*$ is a right $\K$-module by
\[
(a_n\cdots a_0)k=
\sum(a_n\al k_{(n)})\cdots(a_0\al k_{(0)})
\qquad(a_0,\cdots, a_n\in\m^*,\; k\in\K).
\]
Hence it is sufficient to show that 
$\Gr^\perp$ is a right $\K$-submodule of 
$\m^*\otimes \m^*$.
Let 
$\sum b_r\otimes b'_r\in\Gr^\perp$ and $k\in\K$.
For $\sum x_c\otimes y_c\in\Gr$ we have
\begin{align*}
\sum \langle (b_r\al k_{(1)})\otimes (b'_r\al k_{(0)}),x_c\otimes y_c\rangle
=
\sum \langle b_r\otimes b'_r,\ad(k_{(1)})(x_c)\otimes\ad(k_{(0)})(y_c)\rangle,
\end{align*}
and hence we have only to show 
\[
\sum \ad(k_{(1)})(x_c)\otimes\ad(k_{(0)})(y_c)\in\Gr.
\]
We have
\begin{align*}
&\sum (\ad(k_{(1)})(x_c))(\ad(k_{(0)})(y_c))
=\sum
k_{(3)}x_c(S^{-1}k_{(2)})k_{(1)}y_c(S^{-1}k_{(0)})
\\
=&\sum
k_{(1)}x_cy_c(S^{-1}k_{(0)})
=\ad(k)(\sum x_cy_c)=0.
\end{align*}
\end{proof}

We also use the left $\K$-module structure 
\begin{equation}
\K\otimes\wedge\m^*\to\wedge\m^*
\qquad(k\otimes a\mapsto k\ar a)
\end{equation}
of $\wedge\m^*$ given by
\[
k\ar a=a\al Sk
\qquad(k\in \K).
\]
We define a left $\K$-module structure 
\begin{equation}
\label{eq:K-action}
\K\otimes(O_q(G)\otimes\wedge\m^*)\to O_q(G)\otimes\wedge\m^*
\qquad(k\otimes \omega\mapsto k\ar \omega)
\end{equation}
of $O_q(G)\otimes\wedge\m^*$ by
\[
k\ar(\varphi\otimes a)=\sum k_{(0)}\varphi\otimes (k_{(1)}\ar a).
\]

\begin{lemma}
\label{lem:kaf}
For $a\in\wedge^n\m^*$, $\varphi\in O_q(G)$, $k\in \K$ we have
\[
k\ar(a\varphi)=
\sum (k_{(1)}\ar a)(k_{(2)}\varphi(Sk_{(0)})).
\]
\end{lemma}
%%%%%%%%%%%%%%%%%%
\begin{proof}
We show this by induction on $n$.
In the cases $n=0, 1$ the assertion is easily verified by a straightforward calculation.
Assume 
$n>1$.
It is sufficient to show 
\[
k\ar((a\wedge b)\varphi)
=\sum
(k_{(1)}\ar(a\wedge b))(k_{(2)}\varphi (Sk_{(0)}))
\]
for 
$a\in\wedge^{n-1}\Gm_q^*$, $b\in\Gm_q^*$.
By \eqref{eq:Delta1} and \eqref{eq:Delta2} we have
\begin{align*}
&k\ar((a\wedge b)\varphi)
\\=&
\sum k\ar((a\varphi_{(0)})\wedge(b\varphi_{(1)}))
=\sum (k_{(0)}\ar(a\varphi_{(0)}))\wedge (k_{(1)}\ar(b\varphi_{(1)}))
\\
=&\sum
\left((k_{(1)}\ar a)(k_{(2)}\varphi_{(0)}(Sk_{(0)}))\right)
\wedge
\left((k_{(4)}\ar b)(k_{(5)}\varphi_{(1)}(Sk_{(3)}))\right)
\\
=&\sum
\left((k_{(1)}\ar a)(\varphi_{(0)}(Sk_{(0)}))\right)
\wedge
\left((k_{(4)}\ar b)(k_{(5)}\varphi_{(1)}k_{(2)}(Sk_{(3)}))\right)
\\
=&
\sum
\left((k_{(1)}\ar a)(\varphi_{(0)}(Sk_{(0)}))\right)
\wedge
\left((k_{(2)}\ar b)(k_{(3)}\varphi_{(1)})\right)
\\
=&
\sum
\left((k_{(1)}\ar a)\wedge(k_{(2)}\ar b)\right)(k_{(3)}\varphi(Sk_{(0)}))
\\
=&
\sum
\left(k_{(1)}\ar(a\wedge b)\right)(k_{(2)}\varphi(Sk_{(0)})).
\end{align*}
\end{proof}

The right $U_q(\Gg)$-module structure of 
$O_q(G)\otimes\wedge\m^*$
\begin{equation}
\label{eq:U-action}
(\varphi\otimes a)u=
\varphi u\otimes a
\qquad(\varphi\in O_q(G), a\in\wedge\m^*, u\in U_q(\Gg))
\end{equation}
and the left $\K$-module structure \eqref{eq:K-action}
gives an $(O_q(G),O_q(L_J))$-bicomodule structure of 
$O_q(G)\otimes\wedge\m^*$.

For $n\geqq0$ we define a linear map
\begin{equation}
d:O_q(G)\otimes\wedge^n\m^*\to O_q(G)\otimes\wedge^{n+1}\m^*
\end{equation}
by 
\[
d(\varphi\otimes a)=\sum_{\beta\in\Xi_J}e_\beta\varphi\otimes (e^*_\beta\wedge a)
\qquad(\varphi\in O_q(G), a\in\wedge^n\m^*).
\]
%%%%%%%%%%%%%%%%%%
\begin{theorem}
\label{thm:1}
\begin{itemize}
\item[(i)]
The linear endomorphism $d$ of $O_q(G)\otimes\wedge\m^*$ of degree 1 gives a dg-algbra structure of the graded algebra $O_q(G)\otimes\wedge\m^*$.
Namely, we have $d^2=0$ and 
\begin{align}
\label{eq:PP}
d(\omega\wedge \xi)=&d\omega\wedge\xi+(-1)^n \omega\wedge d\xi
\\
\nonumber
&\qquad(\omega\in O_q(G)\otimes\wedge^n\m^*,\; \xi\in
O_q(G)\otimes\wedge\m^*).
\end{align}
\item[(ii)]
The graded algebra $O_q(G)\otimes\wedge\m^*$ is generated by 
$O_q(G)=O_q(G)\otimes\wedge^0\m^*$ and 
$\Image(d:O_q(G)\otimes\wedge^0\m^*
\to O_q(G)\otimes\wedge^1\m^*)$.
Hence $O_q(G)\otimes\wedge\m^*$ is a differential calculus over $O_q(G)$ in the sense of \cite{BeM}.

\item[(iii)] 
The differential calculus $O_q(G)\otimes\wedge\m^*$ is 
$(O_q(G),O_q(L))$-bicovariant.
Namely, the multiplication 
\[
(O_q(G)\otimes\wedge\m^*)\otimes (O_q(G)\otimes\wedge\m^*)\to O_q(G)\otimes\wedge\m^*
\]
and $d:O_q(G)\otimes\wedge\m^*\to O_q(G)\otimes\wedge\m^*$ are homomorphisms of $(O_q(G),O_q(L))$-bicomodules.
\end{itemize}
\end{theorem}
\begin{proof}
(i) By Lemma \ref{lem:eei} we have $d^2=0$.

Let us show \eqref{eq:PP}.
Let us first show
\begin{equation}
\label{eq:m0n0}
d(\varphi\psi)=d\varphi\wedge\psi+\varphi\wedge d\psi
\qquad
(\varphi, \psi\in O_q(G)=O_q(G)\otimes \wedge^0\m^*).
\end{equation}
Identifying $O_q(G)\otimes\m^*$ with $O_q(G)\otimes(O_q(G)^+/\SI_0)$ we have
\begin{equation}
\label{eq:d}
d\varphi=\sum\varphi_{(0)}\otimes[\varphi_{(1)}^+]
\qquad(\varphi\in O_q(G)).
\end{equation}
In fact for $\beta\in\Xi_J$, $u\in U_q(\Gg)$ we have
\begin{align*}
&\sum\langle\varphi_{(0)},u\rangle\langle\varphi_{(1)}^+,e_\beta\rangle=
\sum\langle\varphi_{(0)},u\rangle\langle\varphi_{(1)},e_\beta\rangle
=
\langle \varphi,ue_\beta\rangle=\langle e_\beta\varphi,u\rangle=
\langle d\varphi, u\otimes e_\beta\rangle.
\end{align*}

Therefore we have
\begin{align*}
&d(\varphi\psi)-d\varphi\wedge\psi
\\
=&
\sum\varphi_{(0)}\psi_{(0)}\otimes[(\varphi_{(1)}\psi_{(1)})^+]
-
\sum\varphi_{(0)}\psi_{(0)}\otimes[\varphi_{(1)}^+\psi_{(1)}]
\\
=&
\sum\varphi_{(0)}\psi_{(0)}\otimes
([(\varphi_{(1)}\psi_{(1)}-
\varepsilon(\varphi_{(1)})
\varepsilon(\psi_{(1)})
)-
(\varphi_{(1)}\psi_{(1)}-
\varepsilon(\varphi_{(1)})
\psi_{(1)})
]
\\
=&
\sum\varphi\psi_{(0)}\otimes[\psi_{(1)}^+]
=\varphi\wedge d\psi.
\end{align*}
We see from \eqref{eq:m0n0} that 
\begin{equation}
\label{eq:m0n}
d(\varphi\omega)=d\varphi\wedge\omega+\varphi d\omega
\qquad
(\varphi\in O_q(G), \;
\omega\in O_q(G)\otimes\wedge\m^*).
\end{equation}
In fact, for $\omega=\psi\otimes a$ we have
\[
d(\varphi\omega)=
d(\varphi\psi\otimes a)=d(\varphi\psi)\wedge a
=d\varphi\wedge\psi\wedge a+\varphi d\psi\wedge a
=d\varphi\wedge\omega+\varphi d\omega.
\]
By definition we also have
\begin{equation}
\label{eq:xa}
d(\omega\wedge a)=d\omega\wedge a
\qquad
(\omega\in O_q(G)\otimes\wedge\m^*, 
\;
a\in\wedge\m^*\subset O_q(G)\otimes \wedge\m^*).
\end{equation}
By \eqref{eq:m0n} and \eqref{eq:xa} the proof of 
\eqref{eq:PP} is reduced to showing
\begin{equation}
\label{eq:ap}
d(a\wedge\varphi)=
(-1)^na\wedge d\varphi
\qquad
(a\in\wedge^n\m^*,\; \varphi\in O_q(G)).
\end{equation}
We show \eqref{eq:ap} by induction on $n$.
The case $n=0$ is trivial.
Assume that $n=1$.
We have
\[
d(a\wedge\varphi)=\sum d(\varphi_{(0)}\otimes a\varphi_{(1)})
=\sum e_\beta\varphi_{(0)}\otimes (e^*_\beta\wedge(a\varphi_{(1)})),
\]
\[
a\wedge d\varphi=\sum a\wedge(e_\beta\varphi\otimes e^*_\beta)
=
\sum\varphi_{(0)}\otimes (a(e_\beta\varphi_{(1)})\wedge e^*_\beta).
\]
Hence it is sufficient to show 
\[
\sum \langle e_\beta\varphi_{(0)},u\rangle
e^*_\beta\wedge(a\varphi_{(1)})
=-
\sum\langle\varphi_{(0)},u\rangle
a(e_\beta\varphi_{(1)})\wedge e^*_\beta
\]
for any $u\in U_q(\Gg)$.
This is equivalent to showing
\[
\sum \langle e_\beta\varphi_{(0)},u\rangle
\langle e^*_\beta,x_c\rangle\langle a\varphi_{(1)},y_c\rangle
=-
\sum\langle\varphi_{(0)},u\rangle
\langle a(e_\beta\varphi_{(1)}),x_c\rangle
\langle e^*_\beta,y_c\rangle
\]
for any $u\in U_q(\Gg)$ and $x_c, y_c\in\m$ satisfying $\sum x_cy_c=0$.
This follows from \eqref{eq:Id2} since 
\begin{align*}
&\sum \langle e_\beta\varphi_{(0)},u\rangle
\langle e^*_\beta,x_c\rangle\langle a\varphi_{(1)},y_c\rangle
=
\sum \langle x_c\varphi_{(0)},u\rangle
\langle a\varphi_{(1)},y_c\rangle
\\
=&
\sum \langle \varphi_{(0)},ux_c\rangle
\langle a,y_{c(0)}^+\rangle
\langle \varphi_{(1)},y_{c(1)}\rangle
=
\sum \langle \varphi,ux_cy_{c(1)}\rangle
\langle a,y_{c(0)}^+\rangle,
\\
&\sum\langle\varphi_{(0)},u\rangle
\langle a(e_\beta\varphi_{(1)}),x_c\rangle
\langle e^*_\beta,y_c\rangle
=
\sum\langle\varphi_{(0)},u\rangle
\langle a(y_c\varphi_{(1)}),x_c\rangle
\\
=&
\sum\langle\varphi_{(0)},u\rangle
\langle a,x_{c(0)}^+\rangle
\langle y_c\varphi_{(1)},x_{c(1)}\rangle
=
\sum\langle\varphi,ux_{c(1)}y_c\rangle
\langle a,x_{c(0)}^+\rangle.
\end{align*}
Now assume that $n>1$.
For $a\in\wedge^{n-1}\m^*$, $b\in\m^*$, $\varphi\in O_q(G)$ we have
\begin{align*}
&d(a\wedge b\wedge \varphi)
=\sum d(a\wedge\varphi_{(0)}\wedge(b\varphi_{(1)}))
=\sum d(a\wedge\varphi_{(0)})\wedge(b\varphi_{(1)})
\\
=&(-1)^{n-1}a\wedge d\varphi_{(0)}\wedge(b\varphi_{(1)})
=(-1)^{n-1}a\wedge 
d(\varphi_{(0)}\wedge(b\varphi_{(1)}))
\\
=&
(-1)^{n-1}a\wedge 
d(b\wedge \varphi)
=
(-1)^n a\wedge b\wedge d\varphi
\end{align*}
by \eqref{eq:xa} and the hypothesis of induction.
The proof of \eqref{eq:PP} is complete.

(ii) For $\varphi\in O_q(G)$ we have
\[
\sum(S\varphi_{(0)})\wedge d\varphi_{(1)}
=\sum(S\varphi_{(0)})\varphi_{(1)}\otimes[\varphi_{(2)}^+]
=1\otimes[\varphi^+]
\]
by \eqref{eq:d}, and thus
$O_q(G)\wedge dO_q(G)\supset\m^*$.
Hence the assertion is a consequence of the fact that $O_q(G)\otimes\wedge\m^*$ is generated as an algebra by $O_q(G)$ and $\m^*$.

(iii) It is sufficient to show that 
the multiplication 
\[
(O_q(G)\otimes\wedge\m^*)\otimes (O_q(G)\otimes\wedge\m^*)\to O_q(G)\otimes\wedge\m^*
\]
and $d:O_q(G)\otimes\wedge\m^*\to O_q(G)\otimes\wedge\m^*$ are homomorphisms of $(\K,U_q(\Gg))$-bimodules.
For $\varphi, \psi\in O_q(G)$, $a, b\in\wedge\m^*$, $k\in\K$, $u\in U_q(\Gg)$ we have 
\begin{align*}
&k\ar
((\varphi\otimes a)\wedge(\psi\otimes b))u
\\
=&
k\ar
(\sum\varphi\psi_{(0)}\otimes ((a\psi_{(1)})\wedge b))u
\\
=&
\sum(k_{(0)}\varphi u_{(0)})
(k_{(1)}\psi_{(0)}u_{(1)})\otimes ((k_{(2)}\ar(a\psi_{(1)}))\wedge (k_{(3)}\ar b))
\\
=&
\sum(k_{(0)}\varphi u_{(0)})
(k_{(1)}\psi_{(0)}u_{(1)})\otimes 
(((k_{(3)}\ar a)(k_{(4)}\psi_{(1)}(Sk_{(2)})))\wedge (k_{(5)}\ar b))
\\
=&
\sum(k_{(0)}\varphi u_{(0)})
(\psi_{(0)}u_{(1)})\otimes 
(((k_{(3)}\ar a)(k_{(4)}\psi_{(1)}k_{(1)}(Sk_{(2)})))\wedge (k_{(5)}\ar b))
\\
=&
\sum(k_{(0)}\varphi u_{(0)})
(\psi_{(0)} u_{(1)})\otimes 
(((k_{(1)}\ar a)(k_{(2)}\psi_{(1)}))\wedge (k_{(3)}\ar b))
\\
%=&
%\sum(k_{(0)}\varphi u_{(0)})
%(k_{(2)}\psi)_{(0)}\otimes 
%(k_{(1)}\ar a)(k_{(2)}\psi)_{(1)}\wedge (k_{(3)}\ar b)
%\\
=&
\sum((k_{(0)}\varphi u_{(0)})\otimes(k_{(1)}\ar a))
\wedge
(k_{(2)}\psi u_{(1)}\otimes (k_{(3)}\ar b))
\\
=&
\sum
(k_{(0)}\ar(\varphi\otimes a)u_{(0)})
\wedge
((k_{(1)}\ar(\psi\otimes b)u_{(1)})
\end{align*}
by \eqref{eq:Delta1}, \eqref{eq:Delta2} and Lemma \ref{lem:kaf}, and
\begin{align*}
&k\ar d(\varphi\otimes a)u
=k\ar\sum e_\beta\varphi u\otimes (e^*_\beta\wedge a)
\\
=&\sum k_{(0)}e_\beta\varphi u\otimes
((k_{(1)}\ar e^*_\beta)\wedge(k_{(2)}\ar a))
\\
=&
\sum k_{(0)}(\ad(Sk_{(1)})(e_\beta))\varphi u\otimes
(e^*_\beta\wedge(k_{(2)}\ar a))
\\
=&\sum k_{(0)}(Sk_{(1)})e_\beta k_{(2)}\varphi u\otimes
(e^*_\beta\wedge(k_{(3)}\ar a))
\\
=&
\sum e_\beta k_{(0)}\varphi u\otimes
(e^*_\beta\wedge(k_{(1)}\ar a))
=d(k\ar(\varphi\otimes a)u).
\end{align*}
\end{proof}
%%%%%%%%%%%%%%%%%%
%%%%%%%%%%%%%%%%%%%%%%%
\section{Quantum homogeneous space $O_q(G/L_J)$}
\subsection{}
%%%%%%%%%%%%%%%%%%
We set $O_q(G/L_J)=O_q(G)^{\co(O_q(L_J)}$ for the natural right $O_q(L_J)$-comodule structure of $O_q(G)$
given by the composite of
\[
O_q(G)\xrightarrow{\Delta}O_q(G)\otimes O_q(G)
\xrightarrow{\id\otimes\pi}
O_q(G)\otimes O_q(L_J),
\]
where $\pi:O_q(G)\to O_q(L_J)$ is the canonical Hopf algebra homomorphism.
We have
\begin{align*}
O_q(G/L_J)
=\{f\in O_q(G)\mid kf=\varepsilon(k)f\;(k\in\K)\}.
\end{align*}
We see easily that 
\begin{equation}
\label{eq:Flc}
\Delta(O_q(G/L_J))\subset O_q(G)\otimes O_q(G/L_J),
\end{equation}
so that $O_q(G/L_J)$ is a left coideal subalgebra of $O_q(G)$.
In this section we recall well-known facts related to $O_q(G/L_J)$.

In the rest of this paper we write 
\[
\O=O_q(G),\qquad\bO=O_q(L_J),\qquad
\F=O_q(G/L_J)
\]
for simplicity.

By our assumption on $q$ any surjective homomorphism in
${}^{\bO}\!\Mod$ splits.
Hence we see easily the following (see \cite[Corollary 1.5]{MS})
\begin{lemma}
\label{lem:fc}
The coalgebra $\O$ is faithfully coflat as a right $\bO$-comodule.
\end{lemma}

%%%%%%%%%%%%%%%%%%
\begin{proposition}
\label{prop:QH}
The left coideal subalgebra $\F$ of $\O$ is a quantum homogeneous $\O$-space in the sense of \cite{CO}.
Namely, $\O$ is faithfully flat as a right $\F$-module, 
and 
we have
\begin{equation}
\label{eq:QH1}
\F^+\O=\O\F^+.
\end{equation}
Moreover, we have 
\begin{equation}
\label{eq:QH2}
\bO=\O/\F^+\O.
\end{equation}
\end{proposition}
%%%%%%%%%%%%%%%%%%
\begin{proof}
Masuoka's bijection \cite[Theorem 1.11]{M} together with 
Lemma \ref{lem:fc} implies \eqref{eq:QH2} and
the  faithfully flatness of 
$\O$ as a right $\F$-module.
See \cite[Lemma A.1]{CO} for \eqref{eq:QH1}
\end{proof}
\subsection{}
We denote by ${}^{\O}_{\F}\!\Mod$
the category consisting of left $\O$-comodules equipped with a compatible left $\F$-module structure.
Namely, an object of ${}^{\O}_{\F}\!\Mod$
is a left $\O$-comodule $M$ equipped with a left $\F$-module structure $F\otimes M\to M$
which is a morphism of left $\O$-comodules, 
and a morphism in ${}^{\O}_{\F}\!\Mod$ is a linear map respecting the  left $\O$-comodule structure and  the left $\F$-module structure.

We define functors
\begin{align}
&\Phi:{}^{\O}_{\F}\!\Mod
\to
{}^{\bO}\!\Mod
\\
&\Psi:{}^{\bO}\!\Mod
\to
{}^{\O}_{\F}\!\Mod
\end{align}
by 
\[
\Phi(M)=M/\F^+M, 
\qquad
\Psi(V)=\O\Box_{\bO}V.
\]
Here, the left $\bO$-comodule structure of $\Phi(M)$ is given by
\[
\Phi(M)\to \bO\otimes\Phi(M)
\qquad([m]\mapsto \sum[m_{(-1)}]\otimes [m_{(0)}]),
\]
and the left $\O$-comodule structure and the left $\F$-module structure of $\Psi(V)$ is given by
\[
\Psi(V)\to \O\otimes\Psi(V)
\qquad(
\sum \varphi_r\otimes v_r\mapsto \sum\varphi_{r(0)}\otimes\varphi_{r(1)}\otimes v_r)
\]
and 
\[
\F\otimes\Psi(V)\to\Psi(V)
\qquad
(f\otimes(\sum\varphi_r\otimes v_r)\mapsto\sum f\varphi_r\otimes v_r)
\]
respectively.
Then $\Psi$ is right adjoint to $\Phi$.

By Takeuchi \cite{Tak} together with Proposition \ref{prop:QH} we have the following.
%%%%%%%%%%%%%%%%%%
\begin{proposition}
\label{prop:Tak0}
The functors $\Phi$ and $\Psi$ give the equivalence 
\begin{equation}
\label{eq:Tak0}
{}^{\O}_{\F}\!\Mod
\cong
{}^{\bO}\!\Mod
\end{equation}
of categories.
\end{proposition}
%%%%%%%%%%%%%%%%%%
We denote by ${}^{\O}_{\F}\!\Mod_\F$
the category consisting of left $\O$-comodules equipped with a compatible $\F$-bimodule structure.
We also denote by ${}^{\bO}\Mod_F$ 
the category consisting of left $\bO$-comodules equipped with a compatible right $\F$-module structure.
By Proposition \ref{prop:Tak0} we easily obtain  
the equivalence
\begin{equation}
\label{eq:Tak1}
{}^{\O}_{\F}\!\Mod_\F
\cong
{}^{\bO}\!\Mod_\F
\end{equation}
of categories.

For $M, N\in{}^{\O}_{\F}\!\Mod_\F$ we have
$M\otimes_{\F}N\in{}^{\O}_{\F}\!\Mod_\F$ 
by
\[
f(m\otimes n)f'=(fm)\otimes (nf')
\qquad(m\in M, n\in N, f, f'\in\F),
\]
\[
\sum(m\otimes n)_{(-1)}\otimes (m\otimes n)_{(0)}
=\sum m_{(-1)}n_{(-1)}\otimes m_{(0)}\otimes n_{(0)}
\qquad(m\in M, n\in N).
\]
By this the category ${}^{\O}_{\F}\!\Mod_\F$
turns out to be a monoidal category.

We denote by 
${}^{\bO}\!\Mod_0$
the full subcategory of ${}^{\bO}\!\Mod_\F$ consisting of $V\in{}^{\bO}\!\Mod_\F$ such that $vf=\varepsilon(f)v$ for any $v\in V$ and $f\in\F$.
We have a natural equivalence
\[
{}^{\bO}\!\Mod_0\cong {}^{\bO}\!\Mod
\]
of categories.
Denote by ${}^{\O}_{\F}\!\Mod_0$ the full subcategory of ${}^{\O}_{\F}\!\Mod_\F$ 
corresponding to ${}^{\bO}\!\Mod_0$ under \eqref{eq:Tak1}.
We have the equivalence 
\begin{equation}
\label{eq:Tak2}
{}^{\O}_{\F}\!\Mod_0
\cong
{}^{\bO}\!\Mod
\end{equation}
of categories.
Note that ${}^{\O}_{\F}\!\Mod_0$  is a monoidal subcategory of ${}^{\O}_{\F}\!\Mod_\F$.
Since $\bO$ is a Hopf algebra, the category 
${}^{\bO}\!\Mod$ is also a monoidal category with respect to the usual tensor product of left $\bO$-comodules.
We see easily the following.
%%%%%%%%%%%%%%%%%%
\begin{proposition}
\label{prop:Tak2}
The equivalence \eqref{eq:Tak2} respects the monoidal structures.
For $M, N\in {}^{\O}_{\F}\!\Mod_0$
the isomorphism
\[
\Phi(M\otimes_F N)\cong\Phi(M)\otimes\Phi(N)
\]
is given by
\[
[m\otimes n]\longleftrightarrow[m]\otimes[n]
\]
for $m\in M$, $n\in N$,and for $V, W\in {}^{\bO}\!\Mod$
the isomorphism
\[
\Psi(V)\otimes_F\Psi(W)\cong\Psi(V\otimes W)
\]
is given by
\[
(\sum \varphi_r\otimes v_r)\otimes(\sum \psi_s\otimes w_s)
\longleftrightarrow
\sum \varphi_r\psi_s\otimes(v_r\otimes w_s)
\]
for 
$\varphi_r\otimes v_r\in \O\Box_{\bO}V$, 
$\psi_s\otimes w_s\in \O\Box_{\bO}W$.
\end{proposition}
%%%%%%%%%%%%%%%%%%
%%%%%%%%%%%%%%%%%%%%%%%
\subsection{}

%%%%%%%%%%%%%%%%%%

Define a quotient coalgebra $\bU$ of $U_q(\Gg)$ by
\begin{equation}
\bU=U_q(\Gg)/U_q(\Gg)\K^+.
\end{equation}
The pairing \eqref{eq:dual} for $\Ga=\Gg$ induces a pairing
\begin{equation}
\label{eq:dual2}
\langle\;,\;\rangle:\F\times\bU\to\BC.
\end{equation}
This is non-degenetate by \cite[Proposition 7]{HK2}.

By Lemma \ref{lem:dec} the composite of 
$
\m\hookrightarrow U_q(\Gg)\twoheadrightarrow \bU
$
is injective. 
Since $\m$ is finite-dimensional, the composite of 
$\F\to\overline{U}^*\to\m^*$ induced by \eqref{eq:dual2}  
and the injective linear map $\m\to\bU$ is surjective.
Hence we have
\begin{equation}
\m^*\cong \F^+/\SI
\end{equation}
for a right ideal $\SI$ of $\F$ contained in $\F^+$.

Note that $\bU$ is a left $U_q(\Gg)$-module by
\begin{equation}
\label{eq:KU}
u[x]=[ux]\qquad(u, x\in U_q(\Gg)).
\end{equation}
For $k\in \K$ and $x\in U_q(\Gg)$ we have
\[
[\ad(k)(x)]=k[x]\qquad(k\in\K,\; x\in U_q(\Gg))
\]
in $\bU$.
%%%%%%%%%%%%%%%%%%
%%%%%%%%%%%%%%%%%%
\section{Differential calculus over $O_q(G/L_J)$}
\begin{lemma}
\label{lem:trivial}
The right $\O$-module $\wedge\m^*$ is trivial as a right $\F$-module.
Namely, we have
\[
af=\varepsilon(f)a\qquad
(a\in\wedge\m^*, \; f\in\F).
\]
\end{lemma}
\begin{proof}
By Lemma \ref{lem:mO} (b) and \eqref{eq:Flc} it is sufficient to show that $\m^*$ is trivial as a right $\F$-module.
For $a\in \m^*$, $f\in \F$, $x\in\m$ we have
\[
\langle af,x\rangle
=\sum\langle a,x_{(0)}^+\rangle\langle f, x_{(1)}\rangle
=\sum \langle a,x_{(0)}^+\rangle\varepsilon(f)\varepsilon(x_{(1)})
=\varepsilon(f)\langle a,x\rangle
\]
by Lemma \ref{lem:Dm}.
\end{proof}
Since $\O\otimes\wedge\m^*$ is an $(\O,\bO)$-bicovariant dg-algebra, we obtain a left $\O$-covariant dg-algebra
\[
\Omega=\bigoplus_{n\geqq0}\Omega^n,
\]
where
\begin{align*}
\Omega^n=&(\O\otimes\wedge^n\m^*)^{\co(\bO)}
\\
=&
\{\omega\in \O\otimes\wedge^n\m^*\mid
k\ar\omega=\varepsilon(k)\omega\;(k\in\K)\}.
\end{align*}
Regarding $\wedge\m^*$ as a left $\bO$-comodule via
the right $\K$-module structure of $\wedge\m^*$ given in Lemma \ref{lem:mK} we can also write
\[
\Omega
=
\O\Box_{\bO}(\wedge\m^*).
\]
By Lemma \ref{lem:trivial} we have
\[
\Omega^n\in{}^{\O}_{\F}\!\Mod_0
\]
for any $n\geqq0$.

We have 
\begin{align}
\label{eq:Omega01}
\Omega^0=
\O\Box_{\bO}\BC=\F,
\qquad
\Omega^1=
\O\Box_{\bO}\m^*.
\end{align}

The multiplication 
\[
\mu_{nm}:\Omega^n\otimes_F\Omega^m\to \Omega^{n+m}
\qquad(n, m\geqq0)
\]
of the graded algebra $\Omega$ is a morphism 
in ${}^{\O}_{\F}\!\Mod_0$, and hence 
induces a morphism
\[
\Phi(\mu_{nm}):\Phi(\Omega^n)\otimes\Phi(\Omega^m)\to \Phi(\Omega^{n+m})
\]
in ${}^{\bO}\!\Mod$.
This defines a graded algebra structure of 
$\Phi(\Omega)=\bigoplus_n\Phi(\Omega^n)$.
By $\Omega^n=\O\Box_{\bO}(\wedge^n\m^*)$ 
we have the identification 
$\Phi(\Omega^n)=\wedge^n\m^*$.
Hence $\wedge \m^*$ is endowed with a new algebra structure given by $\Phi(\mu_{nm})$.
The new multiplication of $\wedge \m^*$ is written as
\begin{equation}
\wedge \m^*\times \wedge \m^*\to \wedge \m^*
\qquad((a,b)\mapsto a\triangle b).
\end{equation}
%%%%%%%%%%%%%%%%%
\begin{proposition}
\label{prop:Omega0}
\begin{itemize}
\item[(i)]
The algebra $\Omega$ is generated by $\Omega^0$ and $\Omega^1$.
\item[(ii)]
The multiplication of $\Omega$ is given by
\begin{equation}
\label{eq:Mu}
(\sum\varphi_r\otimes a_r)\wedge
(\sum\psi_s\otimes b_s)=
(-1)^{nm}\sum\varphi_r\psi_s\otimes b_s\wedge a_r
\qquad
\end{equation}
for
\[
\sum\varphi_r\otimes a_r\in\Omega^n,\qquad
\sum\psi_s\otimes b_s\in\Omega^m.
\]
\item[(iii)]
The multiplication of  $(\wedge\m^*,\triangle)$ is given by
\begin{equation}
a\triangle b=(-1)^{nm}b\wedge a
\qquad(a\in\wedge^n\m^*,\; b\in\wedge^m\m^*).
\end{equation}
\end{itemize}
\end{proposition}
%%%%%%%%%%
\begin{proof}
Let
$\Omega_\dagger=\bigoplus_{n\geqq0}\Omega_\dagger^n$ be the graded subalgebra of $\Omega$ generated by $\Omega^0$ and $\Omega^1$.
We first show \eqref{eq:Mu} for 
\begin{equation}
\label{eq:cond}
\sum\varphi_r\otimes a_r\in\Omega_\dagger^n,\qquad
\sum\psi_s\otimes b_s\in\Omega_\dagger^m
\end{equation}
by induction on $m$.
In the case $m=0$ this follows from Lemma \ref{lem:trivial}.
In the case $m=1$ it is sufficient to show 
\begin{equation}
\label{eq:Mu1}
\sum\psi_{s(0)}\otimes (a\psi_{s(1)})\wedge b_s
=
(-1)^{n}\sum\psi_s\otimes b_s\wedge a
\qquad
\end{equation}
for
\[
a\in\wedge^n\m^*,\qquad
\sum\psi_s\otimes b_s\in\O\Box_\bO\m^*.
\]
Let us show \eqref{eq:Mu1} by induction on $n$.
This is obvious for $n=0$.
In the case $n=1$
it is sufficient to show 
\[
\sum\psi_{s(0)}\otimes (a\psi_{s(1)})\otimes b_s
+\sum\psi_s\otimes b_s\otimes a
\in
\O\otimes\Gr^\perp.
\]
For $\sum x_c\otimes y_c\in\Gr$
we have
\begin{align*}
&\sum\langle
(a\psi_{s(0)})\otimes b_s,x_c\otimes y_c\rangle
\psi_{s(1)}
+\sum\langle b_s\otimes a,x_c\otimes y_c\rangle
\psi_s
\\
=&
\sum
\langle a\psi_{s(0)},x_c\rangle
\langle b_s,y_c\rangle
\psi_{s(1)}
+
\sum
\langle b_s,x_c\rangle\langle a,y_c\rangle
\psi_s
\\
=&
\sum
\langle a,x^+_{c(0)}\rangle
\langle \psi_{s(0)},x_{c(1)}\rangle
\langle b_s,y_c\rangle
\psi_{s(1)}
+
\sum
\langle b_s,x_c\rangle\langle a,y_c\rangle
\psi_s
\\
=&
\sum
\langle a,x^+_{c(0)}\rangle
\varepsilon(x_{c(1)}\psi_{s(0)})
\langle b_s,y_c\rangle
\psi_{s(1)}
+
\sum
\langle b_s,x_c\rangle\langle a,y_c\rangle
\psi_s
\\
=&
\sum
\langle a,x^+_{c(0)}\rangle
\varepsilon(\psi_{s(0)})
\langle b_s\al x_{c(1)},y_c\rangle
\psi_{s(1)}
+
\sum
\langle b_s,x_c\rangle\langle a,y_c\rangle
\psi_s
\\
=&
\sum
\langle a,x^+_{c(0)}\rangle
\langle b_s,\ad(x_{c(1)})(y_c)\rangle
\psi_{s}
+
\sum
\langle b_s,x_c\rangle\langle a,y_c\rangle
\psi_s
\\
=&
\sum
\langle b_s\otimes a,
\ad(x_{c(1)})(y_c)\otimes x^+_{c(0)}+
x_c\otimes y_c\rangle
\psi_s
\\
=&0
\end{align*}
by \eqref{eq:Id3}.
Hence \eqref{eq:Mu1} holds  for $n=1$.
Let us show \eqref{eq:Mu1} for $n>1$ assuming the assertion for smaller $n$'s.
Let $a'\in\wedge^{n-1}\m^*$, $a''\in\m^*$.
We have
\[
\sum\psi_{s(0)}\otimes ((a'\wedge a'')\psi_{s(1)})\wedge b_s
=\sum\psi_{s(0)}\otimes((a'\psi_{s(1)})\wedge (a''\psi_{s(2)})\wedge b_s).
\]
By \eqref{eq:Delta1} we see easily that 
\[
\sum\psi_{s(0)}\otimes\psi_{s(1)}\otimes b_s
\in\O\otimes(\O\Box_\bO\m^*).
\]
Thus by \eqref{eq:Mu1} for $n=1$ we  have
\[
\sum\psi_{s(0)}\otimes\psi_{s(1)}\otimes(a''\psi_{s(2)})\wedge b_s=-\sum\psi_{s(0)}\otimes \otimes \psi_{s(1)}\otimes b_s\wedge a''.
\]
Hence by the hypothesis of induction we have
\begin{align*}
&\sum\psi_{s(0)}\otimes ((a'\wedge a'')\psi_{s(1)})\wedge b_s
=-\sum\psi_{s(0)}\otimes(a'\psi_{s(1)})\wedge b_s\wedge a''
\\
=&
(-1)^n\sum\psi\otimes b_s\wedge a'\wedge a''.
\end{align*}
\eqref{eq:Mu} for $m=1$ is proved.
Let us show \eqref{eq:Mu} for $m>1$ 
under the condition \eqref{eq:cond} assuming the assertion for smaller $m$'s.
We may assume 
\[
\sum\psi_s\otimes b_s
=(\sum\psi'_{s'}\otimes b'_{s'})\wedge(\sum\psi''_{s''}\otimes b''_{s''})
\]
for
\[
\sum\psi'_{s'}\otimes b'_{s'}\in\Omega_\dagger^{m-1},
\qquad
\sum\psi''_{s''}\otimes b''_{s''}\in\Omega^{1}.
\]
By \eqref{eq:Mu} for $m=1$ we have
\[
\sum\psi_s\otimes b_s
=(-1)\sum\psi'_{s'}\psi''_{s''}\otimes b''_{s''}\wedge b'_{s'}.
\]
By the hypothesis of induction we have
\begin{align*}
&(\sum\varphi_r\otimes a_r)\wedge
(\sum\psi_s\otimes b_s)=
(\sum\varphi_r\otimes a_r)\wedge
(\sum\psi'_{s'}\otimes b'_{s'})\wedge(\sum\psi''_{s''}\otimes b''_{s''})
\\
=&(-1)^{n(m-1)}
(\sum\varphi_r\psi'_{s'}\otimes b'_{s'}\wedge a_r)
\wedge
(\sum\psi''_{s''}\otimes b''_{s''})
\\
=&(-1)^{n(m-1)+n+m-1}
\sum\varphi_r\psi'_{s'}\psi''_{s''}\otimes b''_{s''}\wedge b'_{s'}\wedge a_r
\\
=&(-1)^{nm}
\sum\varphi_r\psi_s\otimes b_s\wedge a_r.
\end{align*}
We have proved \eqref{eq:Mu} under the condition \eqref{eq:cond}.

Consider the subalgebra $\Phi(\Omega_\dagger)=\bigoplus_n\Phi(\Omega_\dagger^n)$ of $\Phi(\Omega)$.
We identify  $\Phi(\Omega_\dagger)$ with a subalgebra of $\wedge\m^*$ with multiplication $\triangle$.
By \eqref{eq:Mu} under the condition \eqref{eq:cond}
we have 
\begin{equation}
\label{eq:Pm}
a\triangle b=(-1)^{nm}b\wedge a
\qquad(a\in\Phi(\Omega_\dagger^n), \;
b\in\Phi(\Omega_\dagger^m)).
\end{equation}

It remains to show (i).
Since $\Phi(\Omega_\dagger)$ contains $\m^*=\wedge^1\m^*$, we see easily from \eqref{eq:Pm} that 
we have $\Phi(\Omega_\dagger)=\Phi(\Omega)$.
It follows from Proposition \ref{prop:Tak0} that the inclusion $\Omega_\dagger\hookrightarrow\Omega$ is surjective.
We are done.
\end{proof}

\begin{corollary}
\label{cor:Omega}
The algebra $(\wedge\m^*,\triangle)$ is isomorphic to the algebra opposite to 
$(\wedge\m^*,\wedge)$ by the correspondence
\[
a\leftrightarrow (-1)^{n(n-1)/2} a^\circ
\qquad(a\in\wedge^n\m^*),
\]
where 
$a$ on the left is regarded as an element of $(\wedge\m^*,\triangle)$ and $a^\circ$ on the right is the corresponding element of the algebra opposite to $(\wedge\m^*,\wedge)$.
\end{corollary}
%%%%%%%%%%%%%%%%%%
\begin{theorem}
\label{thm:2}
The dg-algebra $\Omega$ is a left $\O$-covariant differential calculus over $\F$.
Moreover, it is the maximal prolongation of the first order differential calculus 
$d:\F\to \Omega^1$.
Namely, for any differential calculus $\Omega'$ which is a prolongation of the first order differential calculus 
$d:\F\to\Omega^1$ there exists uniquely a morphism of differential calculus $\Omega\to\Omega'$ which is identity on the degree zero and one parts.
\end{theorem}
\begin{proof}
By \cite{H} $d:\F\to \Omega^1$
is a first order differential calculus.
Hence $\Omega$ is a differential calculus over $\F$
by Proposition \ref{prop:Omega0} (i).

Let $\Omega_\max=\bigoplus_n\Omega^n_\max$ be the maximal prolongation of the first order differential calculus 
$d:\F\to \Omega^1$.
It remains to show that the canonical surjective homomorphism 
\[
p=\bigoplus_np^n:\Omega_\max
\to \Omega
\]
is an isomorphism.
By Proposition \ref{prop:Tak0} it is sufficient to show that the surjective algebra homomorphism $\Phi(p)=\bigoplus_n\Phi(p^n):\Phi(\Omega_\max)\to\Phi(\Omega)$ is an isomorphism.
This follows from \cite[Proposition 3.6]{HK3} in view of Corollary \ref{cor:Omega}; however, 
we give here a more direct proof using \cite{O}
for the convenience of the readers.
By \cite{O} we have
\[
\Phi(\Omega_\max)=T\m^*/(T\m^*)\CR(T\m^*)
\]
with
\[
\CR=\{\sum[Sf_{(0)}^+]\otimes[f_{(1)}^+]\mid f\in\SI\}
\subset
\m^*\otimes\m^*,
\]
where we identify $\m^*\otimes\m^*$ with
$(\O^+/\SI_0)\otimes(\F^+/\SI)$.
Hence we have a surjective algebra homomorphism
$\Phi(p)$ from $T\m^*/(T\m^*)\CR(T\m^*)$ to $(\wedge\m^*,\triangle)$
which is identity on the degree one part $\m^*$.
In particular, we have $\CR\subset P(\Gr^\perp)$, 
where $P:\m^*\otimes\m^*\to\m^*\otimes\m^*$ is the transposition.
By  Corollary \ref{cor:Omega} it is sufficient to show $\CR=P(\Gr^\perp)$.
Therefore, we have only to show 
$\dim\CR\geqq\dim(\m\otimes\m)/\Gr$

%\[
%T_{\F}\Omega^1=\bigoplus_{n=0}^\infty
%\Omega^1\otimes_\F\cdots\otimes_\F\Omega^1
%\]
%of the $\F$-bimodule $\Omega^1$ by the two-sided ideal 
%generated by 
%the subobject $\CR$ of $\Omega^1\otimes_\F\Omega^1$
%in ${}^{\O}_{\F}\!\Mod_0$, 
%which corresponds under ? to the suboject 
%\[
%\overline{\CR}=\{\sum[Sf_{(0)}^+]\otimes[f_{(1)}^+]\mid f\in\SI\}
%\]

For $x, y\in\m$ and $f\in \SI$ we have
\[
\sum\langle Sf_{(0)}^+,x\rangle\langle f_{(1)}^+,y\rangle
=
\sum\langle Sf_{(0)},x\rangle\langle f_{(1)},y\rangle
=
\sum\langle f_{(0)},Sx\rangle\langle f_{(1)},y\rangle
=\langle f,(Sx)y\rangle.
\]
Considering the weights we see easily by Lemma \ref{lem:dec} that the canonical linear map
\[
\BC1\oplus \Gm_q\oplus(S\m)\m\hookrightarrow \bU
\]
is injective.
Hence 
by the non-degeneracy of \eqref{eq:dual2} we see that 
for any $\xi\in((S\m)\m)^*$ there exists $f\in\SI$ such that $f|_{(S\m)\m}=\xi$.
It follows that $\dim \CR\geqq\dim (S\m)\m$.
It remains to show that $\dim (S\m)\m\geqq\dim\m^2$.
For $x, y\in\m$ we have 
\[
(Sx)y\in-(Sx)(Sy)+U_q(\Gb)U_q(\Gb_J)^+=
-S(yx)+U_q(\Gb)U_q(\Gb_J)^+
\]
by Corollary \ref{cor:Dm}, 
and hence
\[
(S\m)\m+U_q(\Gb)U_q(\Gb_J)^+
=S(\m^2)+U_q(\Gb)U_q(\Gb_J)^+.
\]
On the other hand we have
\[
S(\m^2)\cap U_q(\Gb)U_q(\Gb_J)^+
=
S(\m^2+U_q(\Gb_J)^+U_q(\Gb))
=\{0\}
\]
by Lemma \ref{lem:dec}.
Therefore, we have
\begin{align*}
&\dim (S\m)\m
\\
\geqq&
\dim((S\m)\m+U_q(\Gb)U_q(\Gb_J)^+)/U_q(\Gb)U_q(\Gb_J)^+
\\
\geqq&
\dim(S(\m^2)+U_q(\Gb)U_q(\Gb_J)^+)/U_q(\Gb)U_q(\Gb_J)^+
\\
=&\dim S(\m^2)=
\dim\m^2.
\end{align*}
\end{proof}
\begin{remark}
We can similarly show that the differential calculus $O_q(G)\otimes\wedge\m^*$ treated in Section \ref{sec:O} is the maximal  prolongation of the first order differential calculus 
$d:O_q(G)\to O_q(G)\otimes\m^*$.
We need to use the so-called fundamental theorem 
of two-sided Hopf modules instead of Takeuchi's equivalence.
\end{remark}
%%%%%%%%%%%%%%%%%%
\section{Cohomology groups of vector bundles}
Let $V$ be an integrable right $U_q(\Gp_J)$-module.
We associate to it a left $\O$-covariant left $\F$-module
\[
\CV=\O\Box_{\bO}V.
\]
The left $\O$-comodule structure and the left $\F$-module structure are given by
\[
(\sum\varphi_r\otimes v_r)u=
\sum\varphi_ru\otimes v_r,\qquad
f\sum\varphi_r\otimes v_r=\sum f\varphi_r\otimes v_r
\]
for $\sum\varphi_r\otimes v_r\in \O\Box_{\bO}V$, 
$u\in U_q(\Gg)$, $f\in\F$.
Recall that $\Omega$ is a left $\O$-covariant $\F$-bimodule.
Hence we have a left $\O$-covariant left $\F$-module 
$\Omega\otimes_{\F}\CV$.
By Proposition \ref{prop:Tak2} and Lemma \ref{lem:trivial} we have
\begin{equation}
\label{eq:OV}
\Omega\otimes_{\F}\CV\cong 
\O\Box_{\bO}(\wedge\m^*\otimes V),
\end{equation}
where the left $\bO$-comodule structure of 
$\wedge\m^*\otimes V$ is given by 
\[
(a\otimes v)k=\sum (a\al k_{(0)})\otimes v k_{(1)}
\qquad
(a\in \wedge\m^*, v\in V, k\in\K).
\]
For $\omega=\sum\varphi_r\otimes a_r\in \O\Box_{\bO}(\wedge\m^*)=\Omega$ and 
$\upsilon=\sum_s\psi_s\otimes v_s\in \O\Box_{\bO}V=\CV$
we have
\[
\omega\otimes \upsilon\longleftrightarrow
\sum\varphi_r\psi_t\otimes a_r\otimes v_s
\]
under \eqref{eq:OV}.

We define a map
\[
\Omega\otimes_F(\Omega\otimes_F\CV)
\to\Omega\otimes_F\CV
\qquad(\omega\otimes\theta\mapsto\omega\wedge\theta)
\]
by $\omega\wedge(\omega'\otimes\upsilon)
=(\omega\wedge\omega')\otimes\upsilon$
for $\omega, \omega'\in\Omega$, $\upsilon\in\CV$.
Under the identification
\eqref{eq:OV} 
we have
\begin{equation}
\label{eq:we}
(\sum\varphi_r\otimes a_r)\wedge(\sum\psi_s\otimes b_s\otimes v_s)
=(-1)^{nm}\sum\varphi_r\psi_s\otimes (b_s\wedge a_r)\otimes v_s
\end{equation}
for 
\[
\sum\varphi_r\otimes a_r\in
\O\Box_{\bO}\wedge^n\m^*,
\qquad
\sum\psi_s\otimes b_s\otimes v_s\in
\O\Box_{\bO}(\wedge^m\m^*\otimes V).
\]

We can also rewrite \eqref{eq:OV} as 
\[
\Omega\otimes_{\F}\CV\cong 
(\O\otimes V\otimes(\wedge\m^*))^{\co(\bO)},
\]
where $\O\otimes V\otimes(\wedge\m^*)$ is an 
$(\O,\bO)$-bicomodule 
given by
\[
(\varphi\otimes v\otimes a)u
=\varphi u\otimes v\otimes a,\quad
k(\varphi\otimes v\otimes a)
=\sum k_{(0)}\varphi\otimes v(Sk_{(1)})\otimes (k_{(2)}\ar a)
\]
for
$\varphi\in \O$, $v\in V$, $a\in\wedge\m^*$, $u\in U_q(\Gg)$, $k\in\K$.

Define
\[
\tilde{\nabla}_V:\O\otimes V\otimes \wedge\m^*
\to \O\otimes V\otimes \wedge\m^*
\]
by
\[
\tilde{\nabla}_V(\varphi\otimes v\otimes a)
=
\sum_{\beta\in\Xi_J}
e_{\beta(0)}\varphi\otimes v(Se_{\beta(1)})\otimes (e^*_\beta\wedge a)
\quad
(\varphi\in \O, v\in V, a\in\wedge\m^*).
\]
Since $\tilde{\nabla}_V$ is an endomorphism of an $(\O,\bO)$-bicomodule, it induces an endomorphism
\[
\nabla_V:\Omega\otimes_{\F}\CV\to\Omega\otimes_{\F}\CV
\]
of the left $\O$-comodule $\Omega\otimes_{\F}\CV$.
%%%%%%%%%%%%%%%%%%
\begin{theorem}
\label{thm:3}
\begin{itemize}
\item[(i)] We have $\nabla_V^2=0$.
\item[(ii)]
We have
\[
\nabla_V(\omega\wedge \theta)=d\omega\wedge \theta+(-1)^n\omega\wedge\nabla_V(\theta)
\qquad(
\omega\in\Omega^n, \theta\in \Omega\otimes_{\F}\CV).
\]
\end{itemize}
\end{theorem}
\begin{proof}
(i) By Lemma \ref{lem:eei} we have $\tilde{\nabla}_V^2=0$.
Hence $\nabla_V^2=0$.

(ii) The proof is reduced to the special case where $\theta\in\CV=\Omega_0\otimes_F\CV$.
Hence we will show 
\[
\nabla_V(\omega\wedge \upsilon)=d\omega\wedge \upsilon+(-1)^n\omega\wedge\nabla_V(\upsilon)
\qquad(
\omega\in\Omega^n, \upsilon\in\CV).
\]
Write
\[
\omega=\sum\varphi_r\otimes a_r\in
\O\Box_{\bO}\wedge^n\m^*,
\quad
\upsilon=\sum\psi_s\otimes v_s\in
\O\Box_{\bO}V.
\]
Then we have
\begin{align*}
&\nabla_V(\omega\wedge\upsilon)
\\
=&
\sum
\nabla_V(\varphi_r\psi_s\otimes a_r\otimes v_s)
\\
=&
\sum
e_{\beta(0)}(\varphi_r\psi_s)\otimes 
(e^*_\beta\wedge a_r)\otimes v_s(Se_{\beta(1)})
\\
=&
\sum
(e_{\beta(0)}\varphi_r)(e_{\beta(1)}\psi_s)\otimes 
(e^*_\beta\wedge a_r)\otimes v_s(Se_{\beta(2)})
\\
=&
\sum
(e_{\beta(0)}^+\varphi_r)(e_{\beta(1)}\psi_s)\otimes 
(e^*_\beta\wedge a_r)\otimes v_s(Se_{\beta(2)})
+
\sum
\varphi_r(e_{\beta(0)}\psi_s)\otimes 
(e^*_\beta\wedge a_r)\otimes v_s(Se_{\beta(1)})
\\
=&
\sum
(e_{\beta(0)}^+\varphi_r)\psi_s\otimes 
(e^*_\beta\wedge a_r)\otimes v_se_{\beta(1)}(Se_{\beta(2)})
+
\sum
\varphi_r(e_{\beta(0)}\psi_s)\otimes 
(e^*_\beta\wedge a_r)\otimes v_s(Se_{\beta(1)})
\\
=&
\sum
(e_{\beta}\varphi_r)\psi_s\otimes 
(e^*_\beta\wedge a_r)\otimes v_s
+
\sum
\varphi_r(e_{\beta(0)}\psi_s)\otimes 
(e^*_\beta\wedge a_r)\otimes v_s(Se_{\beta(1)}),
\\
&d\omega\wedge\theta
\\
=&(
\sum e_\beta\varphi_r\otimes (e^*_\beta\wedge a_r)
)
\wedge
(\sum\psi_s\otimes v_s)
=
\sum
(e_\beta\varphi_r)\psi_s\otimes (e^*_\beta\wedge a_r)\otimes v_s,
\\
&\omega\wedge\nabla_V(\theta)
\\
=&(\sum\varphi_r\otimes a_r)\wedge
(\sum e_{\beta(0)}\psi_s\otimes (e^*_\beta\wedge b_s)\otimes v_s(Se_{\beta(1)})
\\
=&
(-1)^{n}
\sum
\varphi_r(e_{\beta(0)}\psi_s)\otimes 
(e^*_\beta\wedge a_r)\otimes v_s(Se_{\beta(1)}).
\end{align*}
We are done.
\end{proof}

For $r\geqq0$ we define a left $\O$-comodule
 $H^r(V)$ as the $r$-th cohomology group of the
complex
\[
\Omega\otimes_{\F}\CV
=
[
0\to \Omega^0\otimes_{\F}\CV\xrightarrow{\nabla_V}\Omega^1\otimes_{\F}\CV\xrightarrow{\nabla_V}\
\cdots].
\]

Now we compare $H^r(V)$ with the derived induction functor of \cite{APW}.
Consider the  left exact functor 
\[
\Ind:{}^{O_q(P_J)}\!\Mod\to {}^{\O}\!\Mod
\]
given by
\[
\Ind(V)=\O\Box_{O_q(P_J)}V.
\]
Then we have
\begin{align*}
\Ind(V)\cong&(\O\otimes V)^{\co(O_q(P_J))}
\\
\cong&
\Hom_{U_q(\Gp_J)}(\BC,\O\otimes V)
\\
\cong&
\Hom_{U_q(\Gm_J)}(\BC,\O\otimes V)^{\co(\bO)},
\end{align*}
where $\O\otimes V$ is regarded as a right $O_q(P_J)$-comodule by
\[
y(\varphi\otimes v)=\sum y_{(0)}\varphi\otimes v(Sy_{(1)})
\qquad (\varphi\in \O, v\in V, y\in U_q(\Gp_J)).
\]
By our assumption on $q$ the category $\Mod^{\bO}$ is semisimple, and hence the derived functor of $\Ind$ is given by
\[
R^r\Ind(V)\cong\Ext^r_{U_q(\Gm_J)}(\BC,\O\otimes V)^{\co(\bO)}.
\]
By a general result on quadratic algebras (see \cite{P}, \cite{PP}) we have a free resolution 
\[
\cdots
\xrightarrow{\deru}
U_q(\Gm_J)\otimes(\wedge^1\m^*)^*
\xrightarrow{\deru}
U_q(\Gm_J)\otimes(\wedge^0\m^*)^*
\to
\BC\to 0
\]
of the trivial $U_q(\Gm_J)$-module $\BC$,
where
\[
\deru:U_q(\Gm_J)\otimes(\wedge^n\m^*)^*
\to
U_q(\Gm_J)\otimes(\wedge^{n-1}\m^*)^*
\]
is given by 
\[
\deru(u\otimes x)=\sum_{\beta\in\Xi_J}ue_\beta\otimes xe^*_\beta
\qquad(u\in U_q(\Gm_J), x\in(\wedge^n\m^*)^*).
\]
Here, the right $\wedge\m^*$-module structure of $(\wedge\m^*)^*$ is given by
\[
\langle xa,b\rangle=\langle x, a\wedge b\rangle
\qquad(x\in (\wedge\m^*)^*, a, b\in \wedge\m^*).
\]
Hence $R^r\Ind(V)$ is given by the $r$-th cohomology group of the complex
\begin{align*}
&\Hom_{U_q(\Gm_J)}(U_q(\Gm_J)\otimes(\wedge\m^*)^*,\O\otimes V)^{\co(\bO)}
\\
\cong&
\Hom_{\BC}((\wedge\m^*)^*,\O\otimes V)^{\co(\bO)}
\\
\cong&
(\O\otimes V\otimes(\wedge\m^*))^{\co(\bO)}.
\end{align*}
It is easily checked that this complex is isomorphic to $\Omega\otimes_{\F}\CV$.
Therefore, we obtain the following.
%%%%%%%%%%%%%%%%%%
\begin{theorem}
%\label{thm:}
For an integrable right $U_q(\Gp_J)$-module $V$ we have $H^r(V)\cong R^r\Ind(V)$.
\end{theorem}
%\section*{Acknowledgment}
\bibliographystyle{unsrt}

\end{document}